\documentclass[a4paper,12pt]{article}
\usepackage{tabularx} % extra features for tabular environment
\usepackage{amsmath,amsthm,verbatim,amssymb,amsfonts,amsopn,amscd}
\usepackage{graphicx} % takes care of graphic including machinery
\usepackage[margin=1in,a4paper]{geometry} % decreases margins

\usepackage[square,numbers]{natbib}
\usepackage[utf8]{inputenc}
\usepackage[T1]{fontenc}
\usepackage[UKenglish]{babel}
\usepackage{enumerate}
\usepackage{paralist} % for asparaenum - numbering as paragraphs in proofs
\usepackage[final]{hyperref} % adds hyper links inside the generated pdf file
\hypersetup{
	colorlinks=true,       % false: boxed links; true: colored links
	linkcolor=blue,        % color of internal links
	citecolor=blue,        % color of links to bibliography
	filecolor=magenta,     % color of file links
	urlcolor=blue         
}

\usepackage{tikz-cd}
\tikzcdset{row sep/normal=2.2em}
\tikzcdset{arrow style=math font}

\theoremstyle{plain}
\newtheorem{theorem}{Theorem}[section]
\newtheorem{corollary}[theorem]{Corollary}
\newtheorem{lemma}[theorem]{Lemma}
\newtheorem{proposition}[theorem]{Proposition}

\theoremstyle{definition}
\newtheorem{definition}[theorem]{Definition}
\newtheorem{remark}[theorem]{Remark}

\DeclareMathOperator{\inv}{inv}

\newcommand{\set}[2]{\{{#1}\,\mid\,{#2}\}}

\DeclareMathOperator{\id}{id}

\newcommand{\e}{\varepsilon}

\DeclareMathOperator{\calA}{\mathcal{A}}

\DeclareMathOperator{\DN}{DN}

\DeclareMathOperator{\frakM}{\mathfrak{M}}

\newcommand{\IN}{\mathbb{N}}

\newcommand{\into}{\hookrightarrow}
\newcommand{\monto}{\rightarrowtail}

\newcommand{\la}{\langle}
\newcommand{\ra}{\rangle}
\title{Second-order infinitesimal groups and affine connections}
\author{\href{mailto:filip.bar.research@gmail.com}{Filip B\'{a}r}}
\date{\itshape In memoriam of Bill Lawvere (09.02.1937 - 23.01.2023)}
\begin{document}
  \maketitle

  \begin{abstract}
    \noindent 
    This paper presents new research in infinitesimal algebra by introducing the concept of an infinitesimal group and exploring its properties and ramifications. The author investigates first- and second-order subgroups of Lie groups and demonstrates the use of the second-order infinitesimal group structure to define a Lie bracket of points intrinsic to the Lie group. This construction allows for the derivation of a second-order Baker-Campbell-Hausdorff formula for the infinitesimal group operation and provides a means to reconstruct the Lie bracket of the Lie algebra of a Lie group. The author also characterises all second-order infinitesimal group structures on KL vector spaces as deformations of vector addition by bilinear maps. The main contribution of the paper is the generalisation of the previously established correspondence between symmetric affine connections and second-order infinitesimally affine structures to manifolds with non-symmetric affine connections via second-order infinitesimal groups.
  \end{abstract}
  
\section{Introduction}
  
  \emph{Synthetic Differential Geometry} (SDG) makes precise the notion of infinitesimals in a way that allows us to study and model deep-rooted intuitions as well as the implied constructions based on infinitesimal arguments that led to the development of differential geometry of manifolds. For example, in \cite[thm.~2.2]{Kock:Levi-Civita} Kock has shown that all formal manifolds carry an affine geometry at the infinitesimal level by admitting affine combinations of families of points that are infinitesimally close. This clarified how one could apply the full algebra of affine combinations directly on the manifold rather than having to rely on charts. Kock has made strong and repeated use of this fact in his subsequent work; most notably in \cite{Kock:Synthetic_Geometry_Manifolds}.

  Aiming to understand and apply infinitesimal models of theories beyond the affine case led the author to develop and study the notion of an \emph{infinitesimal model of an algebraic theory} in \cite{Bar:thesis} (see also \cite{Bar:gluing_of_iMATs}\footnote{Whenever we make a reference to \cite{Bar:gluing_of_iMATs} in this paper, we need to consider everything over a Grothendieck topos instead of the category of sets as done in \cite{Bar:gluing_of_iMATs} to be applicable to SDG. As remarked in \cite{Bar:gluing_of_iMATs} all results generalise to this case by utilising the internal language and logic of a Grothendieck topos. (See \cite{Bar:thesis} for some of the proofs in this more general setting.)}). In this framework Kock's result translates as follows: every formal manifold is a natural infinitesimal model of the algebraic theory of affine combinations; it carries a natural (nil-square) \emph{infinitesimal structure}, which defines the tuples of points that are infinitesimally close, and the affine combinations are defined exactly on these tuples satisfying all the familiar equations of an affine space. 

  Going beyond the nil-square case the author was able to show in \cite{Bar:second_order_affine_structures} that  the models of affine combinations on the second-order infinitesimal structure on a formal manifold correspond to symmetric affine connections. This exemplifies that it is indeed appropriate to speak of a structure when doing infinitesimal algebra, rather than just a property of a formal manifold. 

  However, up to this point the only examples of algebraic theories considered have been affine spaces and vector spaces. In this paper we shall close this gap by studying (non-abelian) infinitesimal groups. We will show that the first- and second-order monads of the neutral element\footnote{The use of the term `monad' in infinitesimal calculus as an `atom' of the continuum goes back to Leibniz, and has been adopted by Kock to refer to the infinitesimal neighbourhood of a point.} are natural infinitesimal subgroups of a Lie group. Moreover, we will be able to recover the basic differential geometric constructs of a Lie group directly from the infinitesimal group structure. 

  For example, in the second-order case we will define a Lie bracket of points as an algebra-like commutator and use it to show that the infinitesimal group law is given by the second-order Baker-Campbell-Hausdorff formula \cite{Hausdorff:symbolische-exponential-formel}. Lifting the Lie bracket of points with the first-order $\log$-$\exp$-bijection to the tangent space at the neutral element we recover the familiar Lie bracket of the Lie algebra of the Lie group. Using the second-order $\log$-$\exp$ bijection to transport the second-order infinitesimal group structure to the Lie algebra leads us to a second class of examples of second-order infinitesimal groups on KL vector spaces, which are truncated formal group laws \cite{Bochner:Formal_Lie_Groups}, \cite{Froehlich:formal_groups}. We shall show that all second-order subgroups of KL vector spaces arise this way. 

  Building on the observation that the Lie bracket of points is the negative of the torsion, we can re-construct the second-order infinitesimal group structure from the left-invariant connection on the Lie group alone. This observation allows us to generalise our main result in \cite{Bar:second_order_affine_structures} to (non-symmetric) affine connections on formal manifolds: we will show that affine connections correspond to infinitesimal group structures on the second-order monad at every point of the manifold. These groups are abelian (and coincide with the second-order infinitesimal linear structure), if and only if the affine connection is symmetric.

%%%%%%%%%%%%%%%%%%%%%%%%%%%%%%%%%%%%%%%%%%%%%%%%%%%%%%%%%%%%%%%%%%%%%%%%%%%%%%%%%%%%%%%%%%%%%%%%%%%%%%

\section{Preliminaries}\label{sec:prelims}

  For the convenience of the reader we will recall the basic notions of infinitesimal algebra of infinitesimally affine spaces and vector spaces in the framework of Synthetic Differential Geometry (SDG) and state the main results relevant for this paper.\footnote{See \cite{Bar:second_order_affine_structures} for a more detailed exposition, and \cite{Bar:gluing_of_iMATs}, \cite{Bar:thesis} for the general theory of infinitesimal models of an algebraic theory. Note that the Kock-Lawvere axiom~(III) and first-order infinitesimal structure do not appear in the previously published work.} As in \cite{Kock:Synthetic_Geometry_Manifolds} and \cite{Bar:second_order_affine_structures} we will follow the na\"{i}ve axiomatic approach to SDG; that is, we will be working over a $\mathbb{Q}$-algebra $R$ satisfying the three subsequent \emph{Kock-Lawvere axioms} for the spaces of first- and second-order infinitesimals $D(n)$ and $D_2(n)$ as well as the infinitesimal space $\tilde{D}_2(2,R^n)$ for all $n\geq 1$: 
  
  \begin{enumerate}[(I)]\label{KL_axioms}
    \item For every map $t:D(n)\to R$ there are unique $a_0,\ldots,a_n\in R$ such that 
    $$
      t(d_1,\ldots,d_n)=a_0+\sum_{j=1}^n a_j d_j,
    $$
    where $D(n)$ denotes the set $D(n)=\set{(d_1,\ldots,d_n)\in R^n}{d_id_j=0,\ 1\leq i,j\leq n}$
    \item  Every map $t:D_2(n)\to R$ is a polynomial function for a uniquely determined polynomial in $R[X_1,\ldots,X_n]$ of total degree $\leq 2$ 
    $$
      t(d_1,\ldots,d_n) = a_0 + \sum_{j=1}^n a_{j} d_j + \sum_{1\leq j\leq k\leq n}a_{jk} d_j d_k, 
    $$
    where $D_2(n)$ denotes the set
    $$
      D_2(n)=\set{(d_1,\ldots,d_n)\in R^n}{\text{any product of three $d_j$ vanishes}}
    $$
    \item\footnote{We will require this axiom only for the characterisation of second-order infinitesimal group structures in theorem~\ref{thm:formal-2nd-order-i-groups}. Note that in well-adapted models of SDG (as well as the Zariski model) the whole Kock-Lawvere axiom scheme is satisfied by $R$, including axioms~(I), (II) and (III) \cite[chap.~1.3]{Kock:Synthetic_Geometry_Manifolds}, \cite{Dubuc:Well-adapted_models}, \cite{Kock:Well-adapted_models}.} Every map $t:\tilde{D}_2(2,R^n)\to R$ is a polynomial function for a uniquely determined polynomial in $R[X_1,\ldots,X_n,Y_1,\ldots,Y_n]$ of total degree $\leq 2$ 
    \begin{align*}
      t(\delta,\e) = & a_0 + \sum_{j=1}^n a_{j} \delta_j + \sum_{1\leq j\leq k\leq n}a_{jk} \delta_j \delta_k\\
      &+ \sum_{j=1}^n b_{j} \e_j + \sum_{1\leq j\leq k\leq n}b_{jk} \e_j \e_k + \sum_{1\leq j\leq k\leq n}c_{jk} \delta_j \e_k\\
    \end{align*}
    where $\tilde{D}_2(2,R^n)$ denotes the set
    $$
    \tilde{D}_2(2,R^n)=\set{(\delta,\epsilon)\in (R^n)^2}{x_{j}x_{k}x_{\ell}=x_{j}x_{k}y_{\ell}=0,\ x,y\in\{\delta,\epsilon\},\ 1\leq j,k,\ell\leq n}
    $$
  \end{enumerate}
  
  A (finite-dimensional) \emph{\textbf{KL vector space}} is an $R$-module $V\cong R^n$.\footnote{All KL vector spaces are asumed to be finite-dimensional in this paper.} We define 
	\begin{align*}
    D(V) &= \set{v\in V}{\phi[v]^2=0\text{ for all bilinear maps }\phi:V^2\to R}, \\
    D_2(V)&=\set{v\in V}{\phi[v]^{3}=0\text{ for all trilinear maps $\phi:V^{3}\to R$}},
  \end{align*}
	and
  \begin{multline*}
    \tilde{D}_2(2,V)=\{(v,w)\in V^2 \mid \phi[x]^3=\phi[x,x,y]=0,\\ x,y\in\{v,w\} \text{ for all trilinear }\phi:V^{3}\to R\},
  \end{multline*}
  where for an $\ell$-linear map $\phi:V^\ell\to W$, $\phi[v]^\ell$ means evaluating $\phi$ on the $\ell$-tuple $(v,\ldots,v)$. In the case of $V=R^n$ we have $D(V)=D(n)$, $D_2(V)=D_2(n)$ \cite[prop.~1.2.2]{Kock:Synthetic_Geometry_Manifolds} and $\tilde{D}_2(2,V)=\tilde{D}_2(2,R^n)$.
  
  An important consequence of axiom~(II) is that every map $f:V\to W$ between KL vector spaces has a unique Taylor representation
	$$ 
    f(P)-f(Q)= \partial f(Q)[P-Q] + \frac{1}{2}\partial^2 f(Q)[P-Q]^2
	$$
	if $P-Q\in D_2(V)$. Here the $R$-linear map $\partial f(P_1)$ denotes the derivative of $f$ at $P_1$ and $\partial^2 f(Q)$ stands for the second derivative of $f$ at $Q$, which is a symmetric $R$-bilinear map $V^2\to W$. In the case of $P-Q\in D(V)$ axiom~(I) yields the first-order Taylor representation 
  $$ 
    f(P)-f(Q)= \partial f(Q)[P-Q].
  $$
  \begin{definition}[i-structure]\label{def:i-structure}
		Let $A$ be a space. An \emph{\textbf{infinitesimal structure (i-structure)}} on $A$ is an $\IN$-indexed family $n\mapsto A\la n\ra\subseteq A^n$ such that
		\begin{enumerate}[(1)]
			\item $A\la 1\ra = A$, $A\la 0\ra=A^0=1$ (the `one point' space, or terminal object)
			\item For every map $h:m\to n$ of finite sets 
			and every $(P_1,\ldots,P_n)\in A\la n\ra$ we have $ (P_{h(1)},\ldots,P_{h(m)})\in A\la m\ra$ 
		\end{enumerate}
	\end{definition}
  An $n$-tuple $(P_1,\ldots,P_n)\in A^n$ that lies in $A\la n\ra$ will be denoted by $\la P_1,\ldots,P_n\ra$ and we shall refer to these points as \emph{infinitesimal neighbours}. A map $f:A\to X$ that preserves i-structure, i.e. $f^n(A\la n\ra)\subseteq X\la n\ra$, is called an \emph{i-morphism}.

  On any KL vector space $V$ the \emph{first neighbourhood of the diagonal} 
	$$
	  \set{(P_1,P_2)}{P_2-P_1\in D(V)} = V\la2\ra
	$$ 
	induces the \emph{\textbf{nil-square i-structure}}
	$$
		V\la m\ra=\set{(P_1,\ldots,P_m)}{(P_i,P_j)\in V\la 2\ra,\ 1\leq i,j \leq m}
	$$    
	This i-structure is \emph{generated} by $V\la 2\ra$ and the largest i-structure containing the first neighbourhood of the diagonal. We shall also define the \emph{\textbf{first-order i-structure}} $V_1$ on $V$ 
  \begin{enumerate}[(1)]
    \item $V_1\la 1\ra = V$, $V_1\la 0\ra = V^0=1$
    \item For $m\geq 2$
    \begin{multline*}
      %V_2\la m \ra = \{(P_1,\ldots, P_m) \in V^m \mid \\ (P_{i_1} - P_{j_1},P_{i_{2}}-P_{j_{2}}, P_{i_{3}}-P_{j_{3}})\in \DN_2(V),\\ \text{ for all } i_\ell, j_\ell \in \{1,\ldots, m\}, 1\leq \ell\leq 3 \}
      V_1\la m \ra = \{(P_1,\ldots, P_m) \in V^m \mid (P_{i_1} - P_{j_1},P_{i_{2}}-P_{j_{2}})\in \DN_1(V),\\ \text{ for all } i_\ell, j_\ell \in \{1,\ldots, m\}, 1\leq \ell\leq 2 \}
    \end{multline*}
  \end{enumerate}
  where $\DN_1(V)$ is the set 
	\begin{multline*}
		\DN_1(V)=\{(v_1,v_2) \in D(V)^{2}\,\mid \text{For any bilinear map }\phi: V^{2} \to R,\ \phi[v_1,v_2]=0\}
    %\DN_1(V)=\{(v_1,v_2) \in D(V)^{2}\,\mid \\ \text{For any trilinear map }\phi: V^{3} \to R,\ \phi[v_1,v_2,v_{3}]=0\}
	\end{multline*}
  The first-order i-structure and the nil-square i-structure agree for $n=2$; i.e. $V\la 2\ra = V_1\la 2\ra$. However, the first-order i-structure is properly contained in the nil-square i-structure, as $\la P,Q,R\ra \in V\la 3\ra$ only implies that $\phi[Q-P,R-P]$ is alternating, but not necessarily $\phi[Q-P,R-P]=0$.\footnote{The first-order i-structure has not appeared explicitly in published work on SDG as of yet. Its relevance will become clear in section~\ref{sec:i-groups}.}  
  
  The \emph{second neighbourhood of the diagonal}
	$$
	\set{(P_1,P_2)}{P_2-P_1\in D_2(V)} = V_2\la 2\ra
	$$
  induces the \emph{\textbf{second-order i-structure}} $V_2$ on $V$
  \begin{enumerate}[(1)]
    \item $V_2\la 1\ra = V$, $V_2\la 0\ra = V^0=1$
    \item For $m\geq 2$
    \begin{multline*}
      %V_2\la m \ra = \{(P_1,\ldots, P_m) \in V^m \mid \\ (P_{i_1} - P_{j_1},P_{i_{2}}-P_{j_{2}}, P_{i_{3}}-P_{j_{3}})\in \DN_2(V),\\ \text{ for all } i_\ell, j_\ell \in \{1,\ldots, m\}, 1\leq \ell\leq 3 \}
      V_2\la m \ra = \{(P_1,\ldots, P_m) \in V^m \mid (P_{i_1} - P_{j_1},P_{i_{2}}-P_{j_{2}}, P_{i_{3}}-P_{j_{3}})\in \DN_2(V),\\ \text{ for all } i_\ell, j_\ell \in \{1,\ldots, m\}, 1\leq \ell\leq 3 \}
    \end{multline*}
  \end{enumerate}
  where $\DN_2(V)$ is the set 
	\begin{multline*}
		\DN_2(V)=\{(v_1,v_2, v_{3}) \in D_2(V)^{3}\,\mid \text{For any trilinear map }\phi: V^{3} \to R,\ \phi[v_1,v_2,v_{3}]=0\}
    %\DN_2(V)=\{(v_1,v_2, v_{3}) \in D_2(V)^{3}\,\mid \\ \text{For any trilinear map }\phi: V^{3} \to R,\ \phi[v_1,v_2,v_{3}]=0\}
	\end{multline*}
  For the i-structures we have the following relationship
  $$
      V_1 \into V,\qquad V_1\into V_2,
  $$
  where all inclusions are strict. Although we have $V\la 2\ra \subset V_2\la 2\ra$, $V\la 3\ra$ is not contained in $V_2\la 3\ra$, if $\dim V>2$. This is why the nil-square i-structure does not embed into the second-order i-structure.
	\begin{proposition}\label{prop:i-morphism}
    Every map $f:V\to W$ between KL vector spaces is an i-morphism for the respective nil-square, first-order and second-order i-structures.
  \end{proposition}
  \begin{proof}
    This follows by direct calculation from the first and second-order Taylor expansions, respectively, and the definition of the first and second-order i-structure \cite[thm.~3.2]{Bar:second_order_affine_structures}.
  \end{proof}
  \begin{definition}[i-affine space]\label{def:i-affine_space}
		Let $A\la-\ra$ be an \emph{i-structure} on $A$. Set $\calA(n) = \set{(\lambda_1, \ldots , \lambda_n)\in R^n}{\sum_{j=1}^{n}\lambda_j=1}$. The space $A$ is said to be an \emph{\textbf{infinitesimally affine space (i-affine space)}} (over $R$), if for every $n\in\IN$ there are operations
		$$\calA(n)\times A\la n\ra\to A,\qquad ((\lambda_1,\ldots,\lambda_n),\la P_1,\ldots,P_n\ra)\mapsto \sum_{j=1}^{n} \lambda_j P_j$$
		satisfying the axioms
		\begin{itemize}
			\item (\emph{Neighbourhood}) Let $\lambda^k\in \calA(n)$, $1\leq k\leq m$. If $\la P_1,\ldots,P_n \ra\in A \la n\ra$  then
			$$ \bigl(\sum_{j=1}^n \lambda^1_j P_j,\ldots, \sum_{j=1}^n \lambda^m_j P_j\bigr)\in A\la m\ra$$
			\item (\emph{Associativity}) Let $\lambda^k\in \calA(n)$, $1\leq k\leq m$, $\mu \in\calA(m)$ and \\ $\la P_1,\ldots,P_n \ra\in A \la n\ra$. We have 
			$$ \sum_{k=1}^m\mu_k\bigl(\sum_{j=1}^n \lambda^k_j P_j\bigr)=\sum_{j=1}^n \bigl(\sum_{k=1}^m\mu_k\lambda^k_j\bigr)P_j$$
			(Note that the left-hand side is well-defined due to the neighbourhood axiom.)
			\item (\emph{Projection}) Let $n\geq 1$ and let $e^n_k\in R^n$ denote the $k$th standard basis vector for $1\leq k\leq n$. For every  $\la P_1,\ldots,P_n \ra\in A \la n\ra$ it holds
			$$\sum_{j=1}^n (e^n_k)_j P_j = P_k$$
			In particular, we have for $n=1$ that $1P=P$, $P\in A$.   
		\end{itemize}
	\end{definition}
  An i-morphism between i-affine spaces is called an \emph{i-affine map}, if it commutes with taking i-affine combinations. 
  
  \begin{proposition}\label{prop:i-affine_spaces}
    \begin{enumerate}[(1)]
      \item Every KL vector space becomes an i-affine space for the nil-square and first-order i-structure, and any map between KL vector spaces is an i-affine map.
      \item Every KL vector space becomes an i-affine space for the second-order i-structure. However, maps between KL vector spaces are not necessarily i-affine.
    \end{enumerate}  
  \end{proposition}
  \begin{proof}
    This follows by direct calculation from taking the first and second-order Taylor expansions, respectively, and the definition of the second-order i-structure \cite[lem.~2.1]{Kock:Levi-Civita}, \cite[thm.~3.3]{Bar:second_order_affine_structures}.  
  \end{proof}

  An \emph{\textbf{i-vector space}} over $R$ is defined in the same vein by replacing each set of affine combinations $\calA(n)$ with the set of $R$-linear combinations $R^n$. Note that for each $P\in A$ an i-affine space $A$ induces an i-vector space structure on every \emph{monad} $\frakM(P)$
  $$
    \frakM(P) = \set{Q\in A}{\la P,Q\ra \in A\la 2\ra}
  $$
  equipped with the induced i-structure
  $$
		\frakM(P)\la n\ra = \set{(Q_1,\ldots, Q_n)\in \frakM(P)^n }{\la P, Q_1,\ldots, Q_n\ra \in A\la n+1\ra}
	$$
  by defining
  $$
    \sum_{j=1}^n \lambda_j Q_j ;= \bigl(1 - \sum_{j=1}^{n} \lambda_j\bigr)P + \sum_{j=1}^{n} \lambda_j Q_j
  $$
  for any $(\lambda_1,\ldots,\lambda_n)\in R^n$ and $\la Q_1,\ldots,Q_n\ra \in\frakM(P)\la n\ra$. The point $P$ serves as the zero vector in $\frakM(P)$.

  A \emph{\textbf{formal manifold}} $M$ is a space that has a cover of \emph{formally open subspaces} $U\monto M$, which are also formally open subspaces of a (finite-dimensional) KL vector space $V$. We call such subspaces the \emph{charts} of $M$.
  $$
  \begin{tikzcd}
    & U \arrow[ld, "\phi"', tail] \arrow[rd, "\iota", tail] &   \\
  V &                                                       & M
  \end{tikzcd}
  $$
  The property of $U\monto X$ being formally open in a space $X$ means that $U$ is stable under all `infinitesimal motions' at each point. In our case we only really require that for each $P\in U$ and any maps $t:D(n)\to X$, $q:D_2(n)\to X$ with $t(0),s(0)\in U$, the maps $t$ and $s$ factor through $U$ for all $n\geq 1$. (See \cite[I.17]{Kock:Synthetic_Geometry_Manifolds} or \cite[def.~3.2.5]{Bar:thesis} for the general definition.)  
  
  The KL vector space $V$ induces the nil-square, first-order and second-order i-structures on any formally open subspace $\iota:U\monto V$ by restriction (= pullback). The embedding $\iota$ becomes an i-morphism that is \emph{infinitesimally closed} (i-closed), i.e. if $\la P_1,\ldots, P_n\ra\in U\la n\ra$ and 
  $\langle \iota(P_1),\ldots, \iota(P_n),Q\rangle\in V\langle n+1\rangle$ then there exists $P_{n+1}\in U$ such that $\langle P_1,\ldots,P_n,P_{n+1}\rangle$ lies in $U\langle n+1\rangle$ and $\iota(P_{n+1})=Q$. This holds true for ll three i-structures. 

  \begin{lemma}\label{lem:formally_open}
    Let $U\monto V$ be a formally open subspace of a KL vector space $V$.
    \begin{enumerate}[(1)]
      \item $U$ is an i-affine subspace for the induced nil-square, first-order and second-order i-structure.
      \item Any map between formally open subspaces of KL vector spaces is an i-morphism. 
      \item Any map between formally open subspaces of KL vector spaces with the nil-square or first-order i-affine structure is an i-affine map.
    \end{enumerate}
  \end{lemma}
  \begin{proof}
    \begin{enumerate}[(i)]
      \item By proposition~\ref{prop:i-affine_spaces}, $V$ is an i-affine space for all these i-structures. Let $\la P_1,\ldots, P_n\ra\in U\la n\ra$. The neighbourhood axiom implies that $\la \iota(P_1),\ldots,\iota(P_n),\sum_{j=1}^{n}\lambda_{j}\iota(P_j)\ra$ for any affine combination $\lambda\in\calA(n)$. Since $\iota$ is i-closed and a mono this defines operations by affine combinations on the i-structure of $U$ making $\iota$ an i-affine map. Moreover, since $\iota:U\monto V$ \emph{reflects i-structure} by construction, i.e. $\la \iota(P_1),\ldots,\iota(P_n)\ra$ implies $\la P_1,\ldots,P_n\ra$, the i-affine combinations satisfy the neighbourhood axiom. The associativity and projection axioms are satified, too, since $\iota$ is a mono. This shows $U$ an i-affine subspace for all three i-structures.
      
      \item Assertions (2) and (3) can be seen as follows. As $\iota$ is i-closed it contains every monad $\frakM(P)$ for each $P\in U$. We have $\frakM(P)\cong D(V)$ for the nil-square and first-order i-structure, and $\frakM_2(P)\cong D_2(V)$ for the second-order i-structure.\footnote{From now on we shall write $\frakM(P)$ for the monad (together with the i-structure) induced by the nil-square i-structure and $\frakM_1(P)$, respectively $\frakM_2(P)$ for the monad (and i-structure) induced by first-order, respectively second-order i-structure on a formal manifold $M$. Note that in \cite{Kock:Synthetic_Geometry_Manifolds} the subindex $1$ is used for the nil-square i-structure not the first-order i-structure. Although $\frakM(P)$ and $\frakM_1(P)$ agree as spaces, the i-structures are different.}\textsuperscript{,}\footnote{Note that $\frakM_2(P)\la 2\ra \cong \tilde{D}_2(2,V)=D_2(V)\la 2\ra$.} Due to the Kock-Lawvere axioms any map $f:\frakM_{(2)}(P)\to W$ has a unique affine, respectively quadratic extension $f:V\to W$ depending on the i-structure. Such maps are i-morphisms by proposition~\ref{prop:i-morphism}. Moreover, in the case of the nil-square and first-order i-structure, $f$ also preserves i-affine combinations by proposition~\ref{def:i-affine_space}~(1). Since formally open subspaces reflect i-structure all of the stated properties also hold true for the maps between formally open subspaces.    
    \end{enumerate}
  \end{proof}
  
  \begin{proposition}\label{prop:formal_mfds}
    Let $M$ be a formal manifold. 
    \begin{enumerate}[(1)]
      \item $M$ carries the nil-square, first-order and the second-order i-structure induced by its charts, and any map between formal manifolds is an i-morphism for each of the the respective i-structures. We have
      $$
        M_1\into M,\qquad M_1\into M_2.
      $$
      \item $M$ carries a natural i-affine structure on its first-order and nil-square i-structure, and any map between formal manifolds is i-affine.
    \end{enumerate}
  \end{proposition}
  \begin{proof}
    \begin{enumerate}[(1)]
      \item All i-structures are constructed as the joint image of the respective i-structures of all the charts of $M$ (a covering family, i.e. an \emph{atlas} of $M$ is sufficent \cite[rem.~3.5]{Bar:second_order_affine_structures}). By lemma~\ref{lem:formally_open}~(2) any map between formal manifolds is an i-morphism for all i-structures \cite[thm.~3.4]{Bar:second_order_affine_structures}.
      \item Note that charts are stable under pullback; the intersection $U\cap W\monto M$ of two charts $\iota:U\monto M$ and $j:W\monto M$ (as subspaces of $M$) is thus a chart.
      $$
      
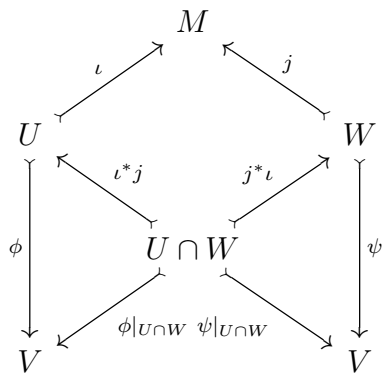
\begin{tikzcd}
        & M                                                                                                                                           &                                                   \\
U \arrow[ru, "\iota", tail] \arrow[dd, "\phi"', tail] &                                                                                                                                             & W \arrow[lu, "j"', tail] \arrow[dd, "\psi", tail] \\
        & U\cap W \arrow[lu, "\iota^*j"', tail] \arrow[ru, "j^*\iota", tail] \arrow[ld, "\phi|_{U\cap W}", tail] \arrow[rd, "\psi|_{U\cap W}"', tail] &                                                   \\
V                                                     &                                                                                                                                             & V                                                
      \end{tikzcd}
      $$
      As can be seen in the diagram $U\cap W$ is the domain of two charts. The identity map $1_{U\cap W}$ induces a map of formally open subspaces $\phi|_{U\cap W}\to \psi|_{U\cap W}$, which is i-affine by lemma~\ref{lem:formally_open}~(3). This shows that the two induced i-affine structures on $U\cap W$ by the two formally open subspaces of $V$ are identical. In addition, the maps $\iota^*j$ and $j^*\iota$ are i-affine and reflect i-structure. The claim now follows from the gluing theorem \cite[thm.~5]{Bar:gluing_of_iMATs}.\footnote{In \cite{Bar:second_order_affine_structures} it is referred to \cite[thm.~3.2.8]{Bar:thesis} to establish (2). However, as pointed out in \cite{Bar:gluing_of_iMATs} the proof of \cite[thm.~3.2.8]{Bar:thesis} relies on \cite[thm.~2.6.19]{Bar:thesis}, which is incorrect.}
    \end{enumerate}  
  \end{proof}
  
  By an \emph{\textbf{affine connection}} on a formal manifold $M$ we shall mean Kock's affine connection on points as defined in \cite[chap.~2.3]{Kock:Synthetic_Geometry_Manifolds} that is based on the geometric idea of completing three points $P,Q,S$ to a parallelogram $PQRS$ with $R=\lambda(P,Q,S)$: An affine connection (on points) is a map $\lambda$ mapping a triple $(P,Q,S)$ with $\la P,Q\ra, \la P,S\ra\in M\la 2\ra$ to a point $\lambda(P,Q,S)$ such that
	\begin{equation}\label{eq:affine_connection}
    \begin{split}
      \lambda(P,Q,P) &= Q\\
	    \lambda(P,P,S) &= S\\
    \end{split}
	\end{equation}
	These properties are sufficient to derive the other nil-square neighbourhood relationships \cite[chap.~2.3]{Kock:Synthetic_Geometry_Manifolds}. An affine connection is called \emph{symmetric}, if 
	$$
    \lambda(P,Q,S)=\lambda(P,S,Q)
  $$
  In a chart $U\monto M$ an affine connection $\lambda$ can be represented as
	$$
		\lambda(P,Q,S)=Q+S-P + \Gamma_P[Q-P,S-P]
	$$
	for a bilinear map $\Gamma_P$ \cite[chapter~2.3]{Kock:Synthetic_Geometry_Manifolds}, which we will refer to as the \emph{connection symbol} of the connection $\lambda$.\footnote{Note that $\Gamma_P$ is the \emph{negative} of the classically defined connection symbol of a covariant derivative. In \cite{Kock:Synthetic_Geometry_Manifolds} the connection symbol is referred to as the \emph{Christoffel symbol}.} The connection is symmetric if and only if its connection symbol is symmetric at every point $P\in M$.

  In SDG the tangent vectors at $P\in M$ are the maps $t:D=D(1)\to M$ with $t(0)=P$. By the Kock-Lawvere axiom~(I) for the space $D$ each tangent vector at $P$ factors through the monad $\frakM(P)$ (induced by the nil-square i-structure) as an i-$R$-linear map. 
  
  \begin{lemma}\label{lem:tanget_vcts_are_i-nbghs}
    For any two tangent vectors $t_1,t_2:D\rightrightarrows M$ at $P$ we have $\la t_1(d), t_2(d)\ra\in M\la 2\ra$ for all $d\in D$.
  \end{lemma}
  \begin{proof}
   \cite[exrc.~2.1.1]{Kock:Synthetic_Geometry_Manifolds} 
  \end{proof}
  
  Due to this lemma the pointwise i-linear structure induced by $\frakM(P)$ on the tangent space at $P$
	$$
		T_PM =\set{t\in M^D}{t(0)=P}
	$$
  makes $T_PM$ a (total) KL vector space \cite[chap.~3.3.2]{Bar:thesis}. By Kock-Lawvere~(I) we get that $T_PM\cong V$, where $V$ is the KL vector space $M$ is modelled on \cite[chap.~4.2]{Kock:Synthetic_Geometry_Manifolds}.

	A symmetric affine connection $\lambda$ induces a \emph{second-order $\log$-$\exp$ bijection}
  $$
    \log_P: \frakM_2(P)\to D_2(T_PM), \qquad
    \exp_P: D_2(T_PM) \to \frakM_2(P),
	$$ 
  which have the following representation in a chart \cite[chap.~8.2]{Kock:Synthetic_Geometry_Manifolds}:
  \begin{equation}\label{eq:log_exp}
    \begin{split}
      \log_P(Q)(d) &= P + d\bigl((Q-P) - \frac{1}{2}\Gamma_P[Q-P]^2\bigr)  \\
      \exp_P(t) &= P + v + \frac{1}{2}\Gamma_P[v]^2
    \end{split}
  \end{equation}
  Note that the tangent vector $t$ is identified with its \emph{principal part} $v\in V$, i.e. the vector $v\in V$ such that $t(d)=P+d\,v$ for all $d\in D$.
  
  \begin{lemma}\label{lem:log_exp}
    \leavevmode \vspace{-\baselineskip}\vspace{1em}
    \begin{enumerate}[(1)]
      \item The first-order $\log_P:\frakM(P)\to D(T_PM)$ and $\exp_P:D(T_PM)\to \frakM(P)$ are i-linear maps.
      \item The second-order $\log_P$ and $\exp_P$ are i-morphisms.
    \end{enumerate}
  \end{lemma}

  \begin{proof}
     Since $T_PM$ is a KL vector space, the first-, respectively, second-order $\log_P$ has an extension $V\to T_PM$ by Kock-Lawvere axioms~(I), respectively (II), when considered in a chart. The assertions follow from propositions~\ref{prop:i-affine_spaces}~(1) and \ref{prop:i-morphism}.
  \end{proof}

  \begin{lemma}\label{lem:embedding-into-2nd-order}
    Let $M$ be a formal manifold and $P\in M$, then
    $
      \frakM(P)\times \frakM(P) \subset \frakM_2(P)\la 2\ra.
    $
  \end{lemma}
  \begin{proof}
    This follows by direct calculation in a chart. (See the proof of \cite[prop.~4.1]{Bar:second_order_affine_structures}.)
  \end{proof}

  \begin{proposition}\label{prop:symmetric_affine_connection}
    Let $M$ be a formal manifold. An i-affine structure on the second-order i-structure on $M$ is equivalent to the existence of a symmetric affine connection $\lambda$ on $M$. 
  \end{proposition}
  \begin{proof}
    Any such i-affine structure defines a symmetric affine connection by setting 
    $$\lambda(P,Q,S) = Q-P+S.$$
    (Note that $P,Q\in\frakM(P)$ implies $\la P,Q\ra\in M_2\la 2\ra$ by the preceding lemma) Conversely, any symmetric affine connection $\lambda$ defines an i-affine structure on $M$ by using the second-order $\log$-$\exp$ bijection (and i-morphisms) to transport the i-linear structure from $D_2(T_PM)$ to $\frakM_2(P)$ for every $P\in M$. 
    $$
      \sum_{j=1}^n\mu_j P_j := \exp_{P}\Bigl(\sum_{j=1}^n\mu_j\,\log_{P}(P_j)\Bigr)
    $$
    for any $\la P_1,\ldots, P_n\ra\in \frakM_2(P)$. For affine combinations this turns out to be independent of the base point $P$. (See the proof of 
    \cite[prop.~4.1, thm.~4.2]{Bar:second_order_affine_structures} for the details.)
  \end{proof}

  Using the representations of the second-order $\log$-$\exp$ bijection~(\ref{eq:log_exp}) induced by $\lambda$ we obtain the following representation of the second-order i-affine combinations in a chart $\iota: U\monto M$ \cite[eq.~(3)]{Bar:second_order_affine_structures}
  \begin{equation}
    \label{eq:2nd-order-iaff-local}
    \begin{split}
      \sum_{j=1}^n\mu_j \iota(P_j) 
        &= \iota\left(\sum_{j=1}^n\mu_j P_j + \frac{1}{2}\Bigl(\Gamma _P\bigl[\sum_{j=1}^n\mu_j\,P_j - P\bigr]^2 - \sum_{j=1}^n\mu_j \Gamma_P[P_j- P]^2\Bigr)\right)			
    \end{split}		
  \end{equation}
  for any $\la P,P_1,\ldots,P_n\ra\in U\la n+1\ra$ and $\mu\in\calA(n)$. Here the i-affine combination on the left hand side is taken in $M$, whereas the i-affine combination on the right hand side is taken in $U$ (applying lemma~\ref{lem:formally_open}~(1)). Note, that the second-order $\log$ and $\exp$ become i-linear maps by construction, extending lemma~\ref{lem:log_exp}.

  The main result of this paper is to extend proposition~\ref{prop:symmetric_affine_connection} to not necessarily symmetric affine connections.

%%%%%%%%%%%%%%%%%%%%%%%%%%%%%%%%%%%%%%%%%%%%%%%%%%%%%%%%%%%%%%%%%%%%%%%%%%%%%%%
\section{Infinitesimal subgroups of Lie groups}\label{sec:i-groups}
 
  So far we have only studied infinitesimal affine spaces and vector spaces. In this section we will define infinitesimal groups and study examples of infinitesimal subgroups of Lie groups. Apart from clarifying and simplifying many familiar basic constructions for Lie groups in differential geometry, the correspondences established for second-order infinitesimal subgroups are key to generalise proposition~\ref{prop:symmetric_affine_connection} to non-symmetric affine connections, which we will prove in section~\ref{sec:affine_connections}.

  \begin{definition}\label{def:i-group}
    Let $G\la-\ra$ be an \emph{i-structure} on $G$. The space $G$ is said to be an \emph{\textbf{infinitesimal group (i-group)}}, if it has two operations $\cdot$ and $\inv$
		$$ \cdot: G\la 2\ra\to G,\quad (P,Q)\mapsto PQ\qquad \inv: G\to G,\quad P\mapsto P^{-1}$$
		as well as a constant $e\in G$ satisfying the axioms
		\begin{itemize}
			\item (\emph{Neighbourhood}) 
        \begin{enumerate}[(I)]
          \item $\la PQ, P_1,\ldots,P_n \ra\in G\la n+1\ra$ for all $\la P,Q, P_1,\ldots,P_n\ra \in G\la n+2\ra$ and $n\geq 0$.
          \item $\la P^{-1}, P_1,\ldots,P_n \ra\in G\la n+1\ra$ for all $\la P, P_1,\ldots,P_n\ra \in G\la n+1\ra$ and $n\geq 0$.
          \item $\la e, P_1,\ldots,P_n \ra\in G\la n+1\ra$ for all $\la P_1,\ldots,P_n\ra \in G\la n\ra$ and $n\geq 0$.
        \end{enumerate} 
      \item (\emph{Associativity}) $(PQ)R = P(QR)$ for all $\la P,Q,R\ra \in G\la 3\ra$.
      \item (\emph{Neutral element}) $e\cdot P = P\cdot e = P$ for all $P\in G$.
			\item (\emph{Inverses}) $P\cdot P^{-1} = P^{-1}\cdot P = e$ for all $P\in G$.
		\end{itemize} 
  \end{definition}
  
  Note that the neighbourhood axiom~(I) guarantees that the iterated products required for associativity are well-defined, and neighbourhood~(III) guarantees that every $P\in G$ can be multiplied with $e\in G$. Due to $\la P,P\ra$ and neighbourhood axiom~(II) we also have $\la P^{-1},P\ra$ for every $P\in G$.

  Recall that the derived formal group operations correspond to formal monomials of the form $X^\alpha = X_1^{\alpha_1}\cdot\ldots\cdot X_n^{\alpha_n}$ with $\alpha\in\mathbb{Z}^n$ and $n\geq 0$. Indeed, the formal monomials of length $n$ are precisely the terms in the term algebra over $n$ variables over the signature of the algebraic theory of groups (see \cite[sec.~3]{Bar:gluing_of_iMATs}, \cite[ch.~1.1.4]{Bar:thesis} and the references cited there).
  Due to the neighbourhood axioms each derived formal group operation has an interpretation in an i-group $G$
  $$
    X^\alpha: G\la n\ra\to G,\quad \la P_1,\ldots,P_n\ra\mapsto P_1^{\alpha_1}\cdot\ldots\cdot P_n^{\alpha_n},
  $$
  where $P^n$ is the $n$-fold product of $P$ with itself, $P^{-n}:=(P^{-1})^n$ and $P^0:=e$ for $n\in\IN$, as in the theory of groups. We shall write $P^\alpha$ for $X^\alpha(P_1,\ldots, P_n)$.
  
  \begin{lemma}\label{lem:iMATs}
    The derived operations $X^\alpha$ satisfy the following neighbourhood axiom:
    Let $\alpha^k\in\mathbb{Z}^n$, $1\leq k\leq m$. If $P=\la P_1,\ldots,P_n \ra\in G \la n\ra$  then
		$ 
       \bigl\langle P^{\alpha^1} ,\ldots, P^{\alpha^m}\bigr\rangle\in G\la m\ra.
    $
    In particular, an i-group $G$ is an infinitesimal model of the algebraic theory of groups as defined in \cite[def.~4]{Bar:gluing_of_iMATs}.
  \end{lemma}
  \begin{proof}
    \begin{enumerate}[(i)]
      \item First note that due to the neighbourhood axioms of an i-group every $X^\alpha$ satisfies $\la P^\alpha, Q_1,\ldots,Q_m \ra\in G\la m+1\ra$ for all $\la P_1,\ldots,P_n,Q_1,\ldots,Q_m\ra \in G\la n+m\ra$ and $m\geq 0$. Indeed, consider the operation $P\mapsto P^k$ for $k>0$ first. Since $\la P ,Q_1,\ldots,Q_m\ra \in G\la 1+m\ra$ we also have 
      $$
      \la \underbrace{P,\ldots,P}_{k\text{ times}}, Q_1,\ldots,Q_m \ra
      $$
      to which we can now iteratively apply the binary i-group operation together with the neighbourhood axiom~(I) $k-1$ times to obtain $\la P^k ,Q_1,\ldots,Q_m\ra \in G\la 1+m\ra$. 
      
      For the general case of $P^\alpha$ we can apply the neighbourhood axioms~(II) and (III) to $\la P_1,\ldots,P_n,Q_1,\ldots,Q_m\ra \in G\la n+m\ra$: First replace each $P_j$ with $e$ if $\alpha_j=0$ and with $P_j^{-1}$ if $\alpha_j<0$; then apply the neighbourhood property for each operation $P_j\mapsto P_j^{|\alpha_j|}$ for each $\alpha_j\neq 0$, $1\leq j\leq n$ iteratively. Finally, apply the binary i-group operation together with the neighbourhood axiom~(I) $n-1$ times. 
    
      \item A similar argument as in part~(i) yields the desired neighbourhood axiom. Firstly, $\la P_1,\ldots,P_n \ra$ implies   
      $$
        \la \underbrace{P_1,\ldots,P_n,\ldots,P_1,\ldots,P_n}_{k\text{ times}} \ra
      $$
      Now apply each operation $X^{\alpha^j}$ for $1\leq j\leq k$ iteratively to obtain $\la P_1^{\alpha^1},\ldots,P_n^{\alpha^k}\ra$ form part~(i). 
    \end{enumerate}
  \end{proof}
  \begin{remark}\label{rem:iMATs}
    In \cite[sec.~4]{Bar:gluing_of_iMATs} we have sketched the \emph{infinitesimalisation of an algebraic theory} as an alternative and more practical approach to obtain infinitesimal models of an algebraic theory given a particular presentation of the theory rather than the approach via clones (see also \cite[def. 2.4.1]{Bar:thesis}). The definition of an i-group provides an example how this can be implemented in practice, and how one can obtain a simplified neighbourhood axiom that guarantees the full neighbourhood axiom schema for all derived operations of the theory. This approach together with lemma~\ref{lem:iMATs} can be easily generalised to any presentation of any algebraic theory.     
  \end{remark}

  Note that if the i-structure on $G$ is generated by $G\la 2\ra$ as in the case of a nil-square i-structure of a formal manifold, for example, then it is sufficient to check the neighbourhood axioms
  \begin{enumerate}[(I)]
    \item $\la PQ, R \ra\in G\la 2\ra$ for all $\la P,Q, R\ra \in G\la 3\ra$,
    \item $\la P^{-1}, R \ra\in G\la 2\ra$ for all $\la P, R\ra \in G\la 2\ra$,
    \item $\la e, R \ra\in G\la 2\ra$ for all $R\in G$.
  \end{enumerate}

  We turn to study examples of i-groups. The first class of examples comes from i-vector spaces, as one would expect: For any formal manifold $M$ and $P\in M$ the monads $\frakM(P)$ and $\frakM_1(P)$ with the induced nil-square, repsectively, first-order i-structure are i-groups for the operations\footnote{Recall that $\frakM(P)$ and $\frakM_1(P)$ agree as spaces, but only differ in the i-structure.}
  \begin{align*}
    +:\frakM_{(1)}(P)\la 2\ra\to\frakM(P),&\qquad (Q,R)\mapsto Q+R-P, \\
    \inv:\frakM(P)\to\frakM(P),&\qquad Q\mapsto 2P-Q  
  \end{align*}
  The neutral element is $P$. Note that in contrast to i-affine spaces i-vector spaces have an implicit constant given by the zero-ary operation $0$. The neighbourhood axiom~(III) therefore follows as a special case of the neigbourhood axiom of an i-vector spaces for the zero-ary operation $0$. In particular, any i-vector space is always a monad of the zero vector.

  If $M$ comes with a symmetric affine connection $\lambda$ then 
  $\frakM_2(P)$ is an i-group for the operations induced by second-order i-affine structure induced by $\lambda$ (cf. proposition~\ref{prop:symmetric_affine_connection}). In particular, we have $Q+R = \lambda(P,Q,R)$ if $Q,R\in\frakM(P)$. This makes $\frakM(P)$ an i-subgroup of $\frakM_2(P)$. We obtain the following inclusions of (abelian) i-groups
  $$
    \frakM_1(P)\into\frakM(P),\qquad \frakM_1(P)\into\frakM_2(P)
  $$
  So far all the examples given were abelian i-groups. The obvious place to look for non-abelian examples are non-abelian \emph{Lie groups} $G$, that is formal manifolds that are also (non-abelian) groups. Note that $G$ is a Lie group if and only if it is an i-group for the \emph{indiscrete i-structure}, i.e. $G\la n\ra=G^n$, $n\in\mathbb{N}$; in other words, if the group operation is defined on all of $G\times G$. The notion of infinitesimal groups subsumes ordinary groups. 

  For any other i-structure on $G$ the neighbourhood axiom~(III) entails that we need to consider the restriction of the group operations to monads of the neutral element $e$. 
  The first natural example to consider would be the monad $\frakM(e)$. Although one can show that $\frakM(e)$ is closed under the group operations, i.e. $\la P,Q\ra$ in $\frakM(e)$ implies $PQ\in\frakM(e)$ and $P\in\frakM(e)$ implies $P^{-1}\in\frakM(e)$ \cite[thm.~6.8.1]{Kock:Synthetic_Geometry_Manifolds}, \emph{the nil-square i-structure fails to satisfy the neighbourhood axioms for $\dim G>2$ and is thus not an i-group, in general.}

  The reason for this failure is that the nil-square i-structure encodes second and higher-order information, which becomes an obstruction to the neighbourhood axiom if and only if the group operation is not abelian. This is why $\frakM(e)$ admits an abelian i-group structure, but not a non-abelian one, in general.
  
  Restricting to $\frakM_1(e)$ yields an i-subgroup of $G$. However, in this case we happen to recover an i-subgroup of the abelian i-group $\frakM(e)$ underlying the i-vector space $\frakM(e)$. This result is unsurprising as it confirms the fact that the group operation of a Lie group is addition in first-order \cite{Bochner:Formal_Lie_Groups}. To obtain non-abelian examples of i-groups we have to consider the second-order monad $\frakM_2(e)$. 

  \newpage 

  \begin{theorem}\label{thm:i-subgroups}
    Let $G$ be a Lie group with neutral element $e\in G$.
    \begin{enumerate}[(1)]
      \item The group operations of $G$ make the monads $\frakM_1(e)$, $\frakM_2(e)$ into infinitesimal subgroups of $G$. We have the subsequent inclusions of i-groups
      $$
        \frakM_1(e)\into\frakM_2(e)\into G
      $$
      \item The i-group structure on $\frakM_1(e)$ coincides with the abelian i-group structure induced by the i-linear structure, i.e. $PQ=P+Q$ for all $\la P,Q\ra\in \frakM_1(e)\la 2\ra$ and $\inv(P)=-P$. In particular, we have an inclusion of abelian i-groups $\frakM_1(e)\into\frakM(e)$.
      \item Let $H$ be a Lie group, then every map $f:G\to H$ that preserves the neutral element induces a map of i-groups $\frakM_1^G(e)\stackrel{f}{\to}\frakM_1^H(e)$.
      \item Let $H$ be a Lie group, then every group homomorphism $f:G\to H$ induces a map of i-groups $\frakM_2^G(e)\stackrel{f}{\to}\frakM_2^H(e)$.
    \end{enumerate}
  \end{theorem}

  The subsequent lemma lies at the heart of the proof of theorem~\ref{thm:i-subgroups}.
  
  \begin{lemma}\label{lem:strong_i-structure}
    Let $\la P_1,Q_1,\ldots, P_n,Q_n\ra\in G_k\la 2n\ra$ then $\la (P_1,Q_1),\ldots, (P_n,Q_n)\ra\in (G\times G)_k\la n\ra$ for all $n\geq 1$ and $k\in\{1,2\}$.
  \end{lemma}
  \begin{proof}(Lemma)
    Let $V$ be the KL vector space $G$ is modelled on. Recall that products of $V$-charts are $V\oplus V$-charts, as formally open subspaces are stable under taking products and $V\oplus V$ is a KL vector space, since $V$ is. In particular, $G\times G$ is a formal manifold \cite[exrc.~2.1.2]{Kock:Synthetic_Geometry_Manifolds}, hence carries all the i-structures by proposition~\ref{prop:formal_mfds}~(1).

    We show the claim for $k=2$. (The case $k=1$ follows from the same argument adapted to $\DN_1$ instead of $\DN_2$.) It suffices to show the assertion for the KL vector space $V$ the Lie group is modelled on, which is done by direct calculation: Let $\phi:(V\times V)^3\to R$ be a trilinear map. Writing
    $$
    \begin{pmatrix} P_{i_k} \\ Q_{i_k} \end{pmatrix}-\begin{pmatrix} P_{j_k} \\ Q_{j_k} \end{pmatrix}
    = \begin{pmatrix} P_{i_k}-P_{j_k} \\ 0 \end{pmatrix} + \begin{pmatrix} 0 \\ Q_{i_k}-Q_{j_k} \end{pmatrix}
    $$
    and expanding
    \begin{equation*}
      \begin{split}
        \phi\left[\begin{pmatrix} P_{i_1} \\ Q_{i_1} \end{pmatrix}-\begin{pmatrix} P_{j_1} \\ Q_{j_1} \end{pmatrix},\begin{pmatrix} P_{i_2} \\ Q_{i_2} \end{pmatrix}-\begin{pmatrix} P_{j_2} \\ Q_{j_2} \end{pmatrix},\begin{pmatrix} P_{i_3} \\ Q_{i_3} \end{pmatrix}-\begin{pmatrix} P_{j_3} \\ Q_{j_3} \end{pmatrix}\right]
      \end{split}  
    \end{equation*}
    for any $i_\ell,j_\ell\in\{1,\ldots,n\}$ and $1\leq\ell\leq 3$ yields a sum of trilinear maps in $P_{i_k}-P_{j_k}$ and $Q_{i_k}-Q_{j_k}$. Since $\la P_1,Q_1,\ldots, P_n,Q_n\ra\in V_2\la 2n\ra$, all these trilinear terms vanish. This shows
    $$
      \left(\begin{pmatrix} P_{i_1} \\ Q_{i_1} \end{pmatrix}-\begin{pmatrix} P_{j_1} \\ Q_{j_1} \end{pmatrix},\begin{pmatrix} P_{i_2} \\ Q_{i_2} \end{pmatrix}-\begin{pmatrix} P_{j_2} \\ Q_{j_2} \end{pmatrix},\begin{pmatrix} P_{i_3} \\ Q_{i_3} \end{pmatrix}-\begin{pmatrix} P_{j_3} \\ Q_{j_3} \end{pmatrix}\right)\in \DN_2(V)
    $$
    and thus $\la (P_1,Q_1),\ldots, (P_n,Q_n)\ra\in (V\times V)_2\la n\ra$ as asserted.

    Since the asserted property holds for $V$, it holds for every formally open subspace $\iota:U\monto V$. This is because $\iota$ preserves and $\iota\times\iota$ reflects i-structure. (The latter is true by definition of the i-structure on $U\times U$.\footnote{Note that the first- and second-order i-structures on $V\times V$, and thus also on $U\times U$, are not the product i-structures.}) But then it holds for every chart and hence for the formal manifold $G$.
  \end{proof}

  %\newpage

  \begin{proof}(Theorem)
    \begin{enumerate}[(i)]
      \item The neighbourhood axiom~(III) of an i-group holds by definition of the monad $\frakM_k(e)$, $k\in\{1,2\}$. It remains to show the neighbourhood axioms~(I) and (II). Let
      $
        \la P,Q,P_1,\ldots,P_n \ra\in\frakM_k(e)\la n+1\ra
      $ 
      then 
      $$
        \la e,e,P,Q,e,P_1,\ldots,e,P_n\ra\in G_k\la 2(n+2)\ra.
      $$
      By lemma~\ref{lem:strong_i-structure} we have $\la (e,e), (P,Q),(e,P_1),\ldots,(e,P_n)\ra$ in $(G\times G)_k$. By proposition~\ref{prop:i-morphism} we find $\la e, PQ,e\cdot P_1,\ldots, e\cdot P_n \ra$ in $G_k$. This shows $\la PQ, P_1,\ldots,P_n\ra\in\frakM_k(e)\la n+1\ra$ and hence neighbourhood axiom~(I). To see $\la \inv(P),P_1,\ldots,P_n\ra$ in $\frakM_k(e)$ apply the i-morphism $m\circ (\inv\times 1_G)$ to $\la (e,P),(P_1,e),\ldots,(P_n,e),(e,e)\ra$ in $(G\times G)_k$, where $m$ denotes the binary group operation. This shows (1).

      \item To see (2) let $\la P,Q\ra$ in $\frakM_1(e)$. We have $\la (P,e),(e,Q),(e,e)\ra$ in $(G\times G)_1$. By proposition~\ref{prop:formal_mfds}~(2) the group operation preserves the i-affine combination $(P,Q)=(P,e)+(e,Q)-(e,e)$, i.e.
      $$
        PQ=P\cdot e + e\cdot Q - e\cdot e = P + Q - e = P+ Q,
      $$
      where the right hand side denotes the vector addition in $\frakM_1(e)$. Similarly, the map $m\circ (\inv\times 1_G)$ is i-affine ($m$ once again denotes the binary group operation) and hence preserves the i-affine combination $(P,P)=(P,e)+(e,P)-(e,e)$, which yields
      $ e = \inv(P)+P-e$
      and thus $\inv(P)=-P$; the right hand side denoting the additive inverse of the vector $P$ in the i-vector space $\frakM_1(e)$.

      \item (3) is an immediate consequence of proposition~\ref{prop:formal_mfds}~(2) and (2); whereas (4) is a direct consequence of proposition~\ref{prop:formal_mfds}~(1). 
    \end{enumerate}
  \end{proof}

  Our next goal is to find a representation of the second-order group operations on $\frakM_2(e)$ extending theorem~\ref{thm:i-subgroups}~(2).

  A Lie group $G$ comes equipped with two \emph{canonical affine connections} $\lambda_l$ and $\lambda_r$ \emph{by left and right translations}, respectively, defined by
  \begin{align*}
    \lambda_l(P,Q,R) &= QP^{-1}R  \\
    \lambda_r(P,Q,R) &= RP^{-1}Q,
  \end{align*}
  where $\la P,Q\ra$, $\la P,R\ra\in G\la 2\ra$. Note that the space of affine connections on a formal manifold $M$ becomes an affine space under the pointwise i-affine operations induced by the nil-square i-affine structure on $M$. Indeed, using a representation in a chart it is easily seen that for any two affine connections $\lambda,\lambda'$ we have $\la \lambda(P,Q,R), \lambda'(P,Q,R)\ra$ in $M$ \cite[exrc.~2.4.3]{Kock:Synthetic_Geometry_Manifolds}. The induced pointwise i-affine structure on the space of affine connections is thus total. In particular, $G$ has a canonical symmetric affine connection given by the \emph{symmetrisation} of $\lambda_r$ and $\lambda_l$
  $$
    \bar{\lambda}(P,Q,R) = \frac{1}{2}\lambda_l(P,Q,R) + \frac{1}{2}\lambda_r(P,Q,R).
  $$
  By proposition~\ref{prop:symmetric_affine_connection} $G$ carries a \emph{canonical second-order i-affine structure} induced by $\bar{\lambda}$. Let $\Gamma$ be the connection symbol of $\lambda_l$ in a chart $U$, i.e.
  $$
    \lambda_l(P,Q,R) = Q+R-P + \Gamma_P[Q-P,R-P],
  $$
  then
  $$
    \lambda_r(P,Q,R) = R+Q-P + \Gamma_P[R-P,Q-P]
  $$
  and
  $$
    \bar{\lambda}(P,Q,R) = Q+R-P + \frac{1}{2}\bigl(\Gamma_P[Q-P,R-P] + \Gamma_P[R-P,Q-P]\bigr).
  $$
  We shall denote the connection symbol of $\bar{\lambda}$ by $\bar{\Gamma}$. Since $\lambda_l(e,P,Q) = PQ$ we get the chart representation of the group operation on $\frakM_2(e)$ for $P,Q\in\frakM(e)$ as
  \begin{equation}\label{eq:local_2nd-order_group_mul}
    PQ = P + Q - e + \Gamma_e[P-e,Q-e].
  \end{equation}
  which can be re-written as
  \begin{equation}\label{eq:local_2nd-order_group_mul2}
    PQ = P + Q - e + \bar{\Gamma}_e[P-e,Q-e] + \frac{1}{2}\bigl(\Gamma_e[P-e,Q-e] - \Gamma_e[Q-e,P-e]\bigr).
  \end{equation}
  \begin{lemma}\label{lem:local_2nd-order_group_mul}
    The representations~(\ref{eq:local_2nd-order_group_mul}) and (\ref{eq:local_2nd-order_group_mul2}) extend uniquely to $\la P,Q \ra\in \frakM_2(e)\la 2\ra$. 
  \end{lemma}
  \begin{proof}
    Note that the ternary operation on $G_2\la 3\ra$ given by $(P,Q,R)\mapsto QP^{-1}R$ is an extension of $\lambda_l$ satisfying the equations of an affine connection~(\ref{eq:affine_connection}). This and the Kock-Lawvere axiom~(III) result in the subsequent presentation in a chart:\footnote{Since the ternary group operation we use for the connection is defined globally we could prove this lemma from axiom~(II), without having to rely on axiom~(III).} (Recall that $\frakM_2(P)\la 2\ra \cong \tilde{D}_2(2,V)$)
    $$
      QP^{-1}R = Q + R - P + B_P[Q-P,R-P],
    $$
    for a unique bilinear map $B_P$. We have 
    $$
      d_1d_2\, B_P[v,w]=B_P[d_1\, v, d_2\, w] =\Gamma_P[d_1\, v, d_2\, w] = d_1d_2\,\Gamma_P[v,w]
    $$
    for all $d_1,d_2\in D$ and hence $B_P=\Gamma_P$ after cancelling the universally quantified d's.\footnote{This well-known \emph{cancellation lemma} in SDG is a direct consequence of the Kock-Lawvere axiom~(I).}

    The same extension applies to $\lambda_r$ and we still have $\la QP^{-1}R, RP^{-1}Q\ra$ in $G_1$, so can form their symmetrisation leading to a (unique) extension of $\bar{\lambda}$ to $G_2\la 3\ra$. This shows that the representations~(\ref{eq:local_2nd-order_group_mul}) and (\ref{eq:local_2nd-order_group_mul2}) extend uniquely to $\la P,Q\ra\in\frakM_2(e)\la 2\ra$ as claimed.
  \end{proof}
  
  Combining (\ref{eq:local_2nd-order_group_mul2}) with the local representation of the second-order i-affine combinations on $\frakM_2(e)$ (cf. equation~(\ref{eq:2nd-order-iaff-local})) we shall show that the  representation of the group operation on $G_2$ is given by
  $$
    PQ = P + Q - e + \frac{1}{2}(PQ - QP)
  $$
  (Note that this expression is a second-order i-affine combination in $G$ instead of a chart.) Using the induced i-linear structure on $\frakM_2(e)$ we define a \emph{Lie bracket} (of points) for $\la P,Q\ra$ in $\frakM_2(e)$ as
  $$
    [P,Q] = PQ - QP\ (= e + PQ - QP)
  $$
  With it we can generalise the presentation in theorem~\ref{thm:i-subgroups}~(2) to the second-order i-group $\frakM_2(e)$ as follows. (See \cite[thm.~6.8.1]{Kock:Synthetic_Geometry_Manifolds} for the case of the nil-square i-structure.)

\newpage
  \begin{theorem}\label{thm:2nd-order-operations}
    Let $G$ be a Lie group. Consider the i-group $\frakM_2(e)$ with its canconical second-order i-linear structure induced by the canonical second-order i-affine structure on $G$. 
    \begin{enumerate}[(1)]
      \item The Lie bracket (of points) is an alternating i-bilinear map on $\frakM_2(e)$; i.e. $[-,P]$ and $[P,-]$ commute with second-order i-linear combinations of $P_j\in\frakM_2(e)$, $1\leq j\leq n$, for any $P\in\frakM_2(e)$ provided $\la P, P_1,\ldots, P_n\ra\in\frakM_2(e)\la n+1\ra$.

      \item In $\frakM_2(e)$ The group operation is given by 
      \begin{equation}\label{eq:2nd-order-group-mul}
        PQ = P+Q + \frac{1}{2}[P,Q]\qquad \la P,Q\ra\in\frakM_2(e)\la 2\ra
      \end{equation}
      and is an i-biaffine map for the second-order i-affine structure.

      \item Taking the inverse of $P$ amounts to a point reflection in $e$
      $$\inv(P)= -P \qquad P\in \frakM_2(e).$$
    \end{enumerate} 
  \end{theorem}
  \begin{proof}
    \begin{enumerate}[(1)]
      \item Firstly, note that $\la P,Q \ra$ in $\frakM_2(e)$ implies $\la PQ, QP \ra$ in $\frakM_2(e)\la 2\ra$ by the neighbourhood axioms. $[P,Q]$ is thus well-defined. It is clearly alternating. To show that it is i-bilinear, we consider $[P,Q]$ for $\la P,Q \ra$ in $\frakM_2(e)$ in a chart.
      
      The presentation of the second-order i-linear combinations on $\frakM_2(e)$  in a chart (cf. equation~(\ref{eq:2nd-order-iaff-local})) yields:
      $$
        [P,Q] = PQ - QP + e + \frac{1}{2}\bigl(\bar{\Gamma}_e[PQ - QP]^2 - \bar{\Gamma}_e[QP - e]^2 + \bar{\Gamma}_e[PQ - e]^2\bigr)
      $$
      Substituting (\ref{eq:local_2nd-order_group_mul}) and bearing in mind that all terms containing $Q-e$ and $P-e$ in higher than second-order vanish, we can simplify this to 
      \begin{equation}\label{eq:lie-bracket-of-points_chart}
        [P,Q] = e + \Gamma_e[P-e,Q-e] - \Gamma_e[Q-e,P-e],
      \end{equation}
      Let $P_j\in\frakM_2(e)$, $1\leq j\leq n$, and $\la P, P_1,\ldots, P_n\ra\in\frakM_2(e)\la n+1\ra$. Then a similar argument by applying (\ref{eq:2nd-order-iaff-local}) to the second-order i-linear combination 
      $\sum_{j=1}^{n}\lambda_j P_j = e+\sum_{j=1}^n\lambda_j(P_j-e)$ and substituting it into (\ref{eq:lie-bracket-of-points_chart}) yields:
      $$
        [P,\sum_{j=1}^n\lambda_j P_j] = \sum_{j=1}^n\lambda_j[P,P_j]
      $$
      As the Lie bracket is alternating it is i-linear in each argument, and thus i-bilinear as claimed.

      \item The second statement can be shown by a direct calculation in a chart as well. Due to the neighbourhood axiom of an i-group $\la P,Q\ra$ in $\frakM_2(e)$ implies $\la P,Q,[P,Q]\ra\in \frakM_2(e)\la 3\ra$. The right hand side of (\ref{eq:2nd-order-group-mul}) is given by the second-order i-affine combination 
      $$
        P + Q + \frac{1}{2}[P,Q] - \frac{3}{2}e 
      $$
      Applying (\ref{eq:2nd-order-iaff-local}) yields:
      \begin{multline*}
        P + Q - e + \frac{1}{2}([P,Q] - e) + \frac{1}{2}\Bigl(\bar{\Gamma}_e\bigl[P - e + Q - e + \frac{1}{2}\bigl([P,Q] - e\bigr)\bigr]^2 \\
        - \bar{\Gamma}_e[P - e]^2 - \bar{\Gamma}_e[Q - e]^2 - \frac{1}{2}\bar{\Gamma}_e\bigl[[P,Q] - e\bigr]^2 \Bigr)
      \end{multline*}

      in a chart. Substituting (\ref{eq:lie-bracket-of-points_chart}) while bearing in mind that all terms containing $Q-e$ and $P-e$ in higher than second-order vanish the expression simplifies to
      $$
        P + Q - e + \frac{1}{2}\bigl(\Gamma_e[P-e,Q-e] - \Gamma_e[Q-e,P-e]\bigr) + \bar{\Gamma}_e[P-e,Q-e],
      $$
      which is the representation of $PQ$ in a chart (\ref{eq:local_2nd-order_group_mul2}).

      To complete the proof of (2) note that taking sums is i-biaffine an the Lie bracket of points is i-bilinear by (1). The group operation is thus a second-order i-biaffine map. 

      \item Apply (1) to $e = PP^{-1}$.
    \end{enumerate}
  \end{proof}
  
  A more careful analysis of the proof of (1) reveals that we not only have $[P,Q]\in\frakM_2(e)$ but even $[P,Q]\in\frakM_1(e)$ (this follows from (\ref{eq:lie-bracket-of-points_chart}), or the i-bilinearity of $[-,-]$.). Moreover, for $\la P, P_1,\ldots, P_n\ra\in\frakM_2(e)\la n+1\ra$ we have $\la [P,P_1],\ldots,[P,P_j]\ra\in \frakM_1(e)\la n\ra$ and the second-order i-linear combination $\sum_{j=1}^n\lambda_j[P,P_j]$ agrees with the first-order one.

  \begin{proposition}\label{prop:lie-bracket-of-points}
    The Lie bracket of points factors through an i-bilinear map\footnote{The Lie bracket $\frakM_2(e)\la 2\ra \to \frakM_1(e)$ becomes an i-morphism if we consider the following i-structure on $\frakM_2(e)$: $\la (P_1,Q_1),\ldots, (P_n,Q_n)\ra$ in $\frakM_2(e)\la 2\ra$ if $\la P_1,Q_1,\ldots, P_n,Q_n\ra$ in $\frakM_2(e)$.}
    $$
      [-,-]:\frakM_2(e)\la 2 \ra \to\frakM_1(e)\into \frakM_2(e)
    $$
  \end{proposition}

  \begin{corollary}\label{cor:lie-bracket-is-commutator}
    The Lie bracket on $\frakM_2(e)$ is given by the infinitesimal group commutator
    $$
      [P,Q] = PQP^{-1}Q^{-1},\qquad \la P,Q\ra\in\frakM_2(e)\la 2\ra
    $$
  \end{corollary}
  \begin{proof}
    Infinitesimal algebra with theorem~\ref{thm:2nd-order-operations} yields
    \begin{align*}
      PQP^{-1}Q^{-1} &= PQ(QP)^{-1}\\
      &= PQ(-QP) \\ 
      &= PQ(2e-QP)\\
      &= 2PQ - PQQP\\
      &= 2PQ - (PQ + QP + \frac{1}{2}[PQ,QP])\\
      &= PQ - QP + \frac{1}{2}([P,[P,Q]]+[Q,[P,Q]])
    \end{align*}
    
    Since both brackets are i-trilinear on $\frakM_2(e)$ they are constant zero ($e$) maps. This shows
    $$
      PQP^{-1}Q^{-1} = PQ - QP + e = [P,Q]
    $$
    as claimed. 
  \end{proof}
  
  Since i-trilinear maps on $\frakM_2(e)$ are zero-maps, the Lie bracket of points trivially satisfies the Jacobi identity. We obtain that $(\frakM_2(e),[-,-])$ forms an \emph{infinitesimal model of a Lie algebra} in the sense of \cite[def.~4]{Bar:gluing_of_iMATs}, which is equivalent to saying that $\frakM_2(e)$ is an i-vector space together with an i-bilinear map $[-,-]$ satisfying the Jacobi identity and a neighbourhood axiom $\la [P,Q],P_1,\ldots,P_n\ra$ for $\la P,Q,P_1,\ldots, P_n\ra$ in $\frakM_2(e)$; the latter being a consequence of $\frakM_2(e)$ being both, an i-vector space and an i-group.
  
  When restricting the Lie bracket to $\frakM(e)$, proposition~\ref{prop:lie-bracket-of-points} together with lemma~\ref{lem:embedding-into-2nd-order} yield a (total) Lie bracket of points on $\frakM(e)$, for which the Jacobi identity is not trivial anymore.  

  \begin{corollary}\label{cor:lie-bracket-of-points-on-nil-square-monad}
    The second-order i-Lie algebra $(\frakM(e),[-,-])$ induces a (total) Lie bracket of points on $\frakM(e)$
    $$
      [-,-]:\frakM(e)\times\frakM(e)\to \frakM(e)
    $$
    that is i-bilinear and satisfies the Jacobi identity. 
  \end{corollary}
  \begin{proof}
    The Jacobi identity follows from corollary~\ref{cor:lie-bracket-is-commutator} and the Ph. Hall identity for the group commutator (see \cite[prop.~6.8.5f]{Kock:Synthetic_Geometry_Manifolds}).

    The i-bilinearity of the Lie bracket follows from proposition~\ref{prop:lie-bracket-of-points}, but can also be seen directly: Let $P\in\frakM(e)$ and consider the i-affine map $[P,-]$. Since $[P,e]=e$ this map is i-linear for the nil-square i-linear structure on $\frakM(e)$. Moreover, since the Lie bracket is alternating, it is i-bilinear. 
    
    (The latter proof does not rely on $P$ being a nil-square neighbour of any of the $P_j$ we form the i-linear combination of. It is thus a slightly stronger notion of i-bilinearity than we used above.)
  \end{proof}  

  Note that $\frakM(e)$ is not an i-Lie algebra as the neighbourhood axiom fails to hold. (This is for the same reason $\frakM(e)$ fails to be an i-group.) The only i-Lie algebras we have on monads are on $\frakM_2(e)$ and $\frakM_1(e)$; the latter being trivial.
  
  However, the Lie bracket of points can be lifted from $\frakM(e)$ to the total vector space $T_eG$. Unlike with the i-linear structure we cannot do this pointwise, as $d\mapsto [t_1(d),t_2(d)]$ is quadratic in $d\in D$ and thus the zero map for $t_1,t_2\in T_eG$. (Note that $[t_1(d),t_2(d)]$ is well-defined due to lemma~\ref{lem:tanget_vcts_are_i-nbghs}.) Instead we have to lift the Lie bracket of points via the (first-order) $\exp$-$\log$ bijection.
  
  $$
  \begin{tikzcd}
    D(T_eG)\times D(T_eG) \arrow[d, "\exp_e\times\exp_e"'] \arrow[rr, "{\overline{[-,-]}}"] &  & D(T_eG)                        \\
    \frakM(e)\times\frakM(e) \arrow[rr, "{[-,-]}"]                    &  & \frakM(e) \arrow[u, "\log_e"']
  \end{tikzcd}
  $$
  
  \begin{theorem}\label{thm:lie-algebra-from-i-algebra}
    Let $G$ be a Lie group. The (first-order) $\exp$-$\log$ lift of the Lie bracket of points on $\frakM(e)$ induces a total Lie algebra on $T_eG$, which agrees with the Lie bracket of left-invariant vector fields.
  \end{theorem}
  \begin{proof}
    \begin{enumerate}[(i)]
      \item Since both $\exp_e$ and $\log_e$ are i-linear (lemma~\ref{lem:log_exp}~(1)) the lift $\overline{[-,-]}$ of the Lie bracket of points yields an i-bilinear, alternating map on $T_eG$ satisfying the Jacobi identity. In particular, we have 
      $\overline{[t,0]}=\overline{[0,t]}=0$ for $t\in D(T_eG)$. As $T_eG$ is a KL vector space $\overline{[-,-]}$ has a unique extension to a bilinear map
      $$
        [-,-]:T_eG\times T_eG\to T_eG
      $$
      Let $t_1$ and $t_2$ be tangent vectors at $e$. We have
      \begin{align*}
        d_1d_2\, [t_1,t_2] &= [d_1\,t_1,d_2\,t_2] \\
        &= \overline{[d_1\,t_1,d_2\, t_2]} \\
        &= -\overline{[d_2\, t_2,d_1\,t_1]} \\
        &= d_2d_1\, (-[t_2,t_1])
      \end{align*}
      for $d_1,d_2\in D$. As $T_eG$ is KL, we can cancel the $d$'s and obtain $[t_1,t_2]=-[t_2,t_1]$. The Jacobi identity follows from a similar argument. This shows $(T_G,[-,-])$ a Lie algebra.
      
      \item Let $t_1,t_2\in T_eG$. In a chart we find
      $$
        [t_1,t_2](d) = e + d\,(\Gamma_e[v_1,v_2]-\Gamma_e[v_2,v_1]),
      $$
      where $v_j$ is the principal part of $t_j$. Let $d_1,d_2\in D$. By lemma~\ref{lem:embedding-into-2nd-order} we have $\la t_1(d_1),t_2(d_2)\ra$ in $\frakM_2(e)$ and thus
      \begin{align*}
        [t_1,t_2](d_1d_2) &= e + \Gamma_e[d_1\,v_1,d_2\,v_2]-\Gamma_e[d_2\,v_2,d_1\,v_i]\\
        &=e + \Gamma_e[t_1(d_1)-e, t_2(d_2)-e]- \Gamma_e[t_2(d_2)-e, t_1(d_1)-e]\\
        &= [t_1(d_1),t_2(d_2)] \\
        &= t_1(d_1)\,t_2(d_2)\,\bigl(t_2(d_2)\,t_1(d_1)\bigr)^{-1},
      \end{align*}
      where we have applied (\ref{eq:lie-bracket-of-points_chart}) and corollary~\ref{cor:lie-bracket-is-commutator}. 
      By \cite[eqs.~4.9.3, 6.6.1]{Kock:Synthetic_Geometry_Manifolds} the Lie bracket of $t_1$ and $t_2$ considered as left-invariant vector fields agrees with the group commutator; hence both constructions yield the same Lie bracket on $T_eG$. 
      
      Indeed, recall that that the Kock-Lawvere axiom~(I) implies that for a KL vector space $V$ and any map $m:D\times D\to V$ with $m(0,d)=m(d,0)=0$ factors uniquely through the product map $D\times D\to D$ as $t:D\to V$, such that $t(d_1d_2)=m(d_1,d_2)$. This is why it is sufficient to compare both Lie brackets on the product $d_1d_2$ for $d_1,d_2\in D$.
    \end{enumerate}
  \end{proof}

  We can use the second-order $\exp$-$\log$ bijection (\ref{eq:log_exp}) on $G$ to transport the second-order i-group structure from $\frakM_2(e)$ to $D_2(T_eG)$ making the $\exp$-$\log$ bijection an i-group isomorphism. The resulting i-group on $T_eG$ can be presented explicitly:
  
  \begin{corollary}\label{cor:i-group-on-tangent-space}
    Let $G$ be a Lie group. The subsequent operations on $T_eG$
    \begin{align*}
      (t_1,t_2)&\mapsto t_1+t_2+\frac{1}{2}[t_1,t_2] \\
      t&\mapsto -t
    \end{align*}
    and $0\in T_eG$ make $D_2(T_eG)$ (equipped with the induced second-order i-structure) into an i-group isomorphic to $\frakM_2(e)$. 
  \end{corollary}
  \begin{proof}
    Recall that the second-order $\exp$-$\log$ bijection is also an i-linear isomorphism by proposition~\ref{prop:symmetric_affine_connection}. With our construction of the Lie bracket on $T_eG$ in theorem~\ref{thm:lie-algebra-from-i-algebra} and due to theorem~\ref{thm:2nd-order-operations} we find
    \begin{align*}
      \log_e(\exp_e(t_1)\exp_e(t_2)) &= \log_e(\exp_e(t_1)+\exp(t_2)+\frac{1}{2}[\exp_e(t_1),\exp_e(t_2)] \\
      &= t_1+t_2+\frac{1}{2}[t_1,t_2]
    \end{align*}
    and
    $$
      \log_e(\exp_e(t_1)^{-1}) = \log_e(-\exp_e(t_1)) = -t_1.
    $$
    The transported i-group structure by the second-order $\exp$-$\log$ bijection thus yields the asserted operations. In particular $D_2(T_eG)$ satisfies the axioms of an i-group for these operations and the $\exp$-$\log$ bijection induces an i-group isomorphism $D_2(T_eG)\cong\frakM_2(e)$.
  \end{proof}
   
  The second-order i-group structure on $D_2(T_eG)$ stated above is the well-known second-order Baker-Campbell-Hausdorff formula \cite{Hausdorff:symbolische-exponential-formel} and an example of a truncated formal group law \cite{Bochner:Formal_Lie_Groups}. Indeed, for a second-order i-group on a KL vector space $V$ we do not require a full formal group law; any bilinear map $B:V\times V\to V$ will do. This provides us with a second class of examples of i-groups besides the i-subgroups of Lie groups studied above. In fact, it turns out that all the i-group structures on $D_2(V)$ are obtained as deformations of the canonical abelian i-group structure by bilinear maps.  

  \begin{theorem}\label{thm:formal-2nd-order-i-groups}
    Let $V$ be a KL vector space and $B:V\times V\to V$ a bilinear map. Then the operations 
    \begin{align*}
      (v,w)&\mapsto v + w + B[v,w] \\
      v&\mapsto -v + B[v]^2
    \end{align*}
    together with $0\in V$ make $D_2(V)$ into an i-group. This i-group is abelian if and only if $B$ is symmetric. Conversely, any i-group structure on $D_2(V)$ is of this form for a uniquely determined bilinear map $B$.
  \end{theorem}
  \begin{proof}
    \begin{enumerate}[(i)]
      \item The associativity, neutral element, inverse and neighbourhood axioms are easily verified by direct calculation bearing in mind that any third-order terms vanish for $\la u,v,w\ra$ in $D_2(V)$.

      \item To see that the neutral element $e\in D_2(V)$ of the i-group $D_2(V)$ is the zero vector $0$ we coordinatize the situation $V\cong R^n$ and utilise the pointwise $R$-algebra structure on $R^n$. Since $e$ must be a second-oder neighbour of every $\delta\in D_2(V)$, the map 
      $$
        D_2(V)\to V, \qquad \delta\mapsto (\delta- e)^3=-3\delta^2e+3\delta e^2
      $$ 
      is the zero map. The Kock-Lawvere axiom~(II) thus implies $e^2 = e=0$.
      
      \item Recall that $D_2(V)\la 2\ra = \tilde{D}_2(2,V)$, so any i-group operation is a map $m:\tilde{D}_2(2,V)\to V$. By the Kock-Lawvere axiom~(III) there is a unique vector $a_0\in V$, unique linear maps $A_1$ and $B_1$ as well as bilinear maps $A_2$, $B_2$ and $C_2$, such that
      $$
        m(\delta,\e) = a_0 + A_1[\delta] + B_1[\e] + A_2[\delta]^2 + B_2[\e]^2 + C_2[\delta,\e].
      $$
      Due to $m(0,0)=0$, $m(0,\e)=\e$, $m(\delta,0)=\delta$ we find $a_0=0$, $A_1=B_1=\id_V$ and $A_2=B_2=0$. The i-group operation is thus of the desired form for a unique bilinear map $B=C_2$:
      $$
        m(\delta,\e) = \delta + \e  + C_2[\delta,\e].
      $$
      Since $-\delta + C_2[\delta]^2$ is an $m$-inverse for $\delta\in D_2(V)$, the uniqueness of an inverse in an i-group shows 
      $$
        \delta^{-1} = -\delta + C_2[\delta]^2
      $$
      as asserted.
    \end{enumerate}  
  \end{proof}

  Note that the Jacobi identity is not required for the bilinear map $B$ to show the induced operations forming a second-order i-group in theorem~\ref{thm:formal-2nd-order-i-groups}, as it is a third-order relationship. 
  
%%%%%%%%%%%%%%%%%%%%%%%%%%%%%%%%%%%%%%%%%%%%%%%%%%%%%%%%%%%%%%%%%%%%%%%%%%%%%%%%%%%%%%%%%%%%%%%%%%%%%%%%%%%%%%%%%%%%%%%%%%%%% 
\section{Affine connections on manifolds and second-order i-group structures}\label{sec:affine_connections}

  The aim of this section is to generalise \cite[thm.~4.2]{Bar:second_order_affine_structures} (stated as proposition~\ref{prop:symmetric_affine_connection} in this paper) to non-symmetric affine connections by constructing a converese to theorem~\ref{thm:formal-2nd-order-i-groups} in the previous section~\ref{sec:i-groups}. 

  There we started with a Lie group $G$ and obtained a non-symmetric affine connection  $\lambda_l$ from the group structure. The second-order i-affine structure induced by its symmetrisation together with the Lie bracket of points allowed us to represent the group operations up to second-order with a second-order BCH formula. However, the Lie bracket of points is the negative of the torsion of the affine connection $\lambda_l$, so the whole second-order i-group structure on $\frakM_2(e)$ is determined by $\lambda_l$ alone.
  
  Indeed, recall that in SDG the \emph{torsion} $\tau$ of $\lambda_l$ is defined as 
  $$
   \tau_P(Q,R)=\lambda_l(\lambda_l(P,Q,R),Q,R)
  $$
  Using the definition of $\lambda_l$ in terms of the group operations we find
  $$
    \tau_P(Q,R)=Q(QP^{-1}R)^{-1}R = QR^{-1}PQ^{-1}R
  $$
  In particular, we get
  $$
  \tau_e(P,Q^{-1}) = PQP^{-1}Q^{-1} = [P,Q],
  $$
  or, equivalently, using the i-bilinearity of the Lie bracket
  $$
    \tau_e(P,Q) = [P,Q^{-1}] = [P,-Q] = -[P,Q]
  $$
  for $P,Q\in\frakM(e)$.

  We now wish to do the converse, i.e. to construct a second-order i-group structure starting with an affine connection. We state the \emph{main result} of this paper: 

  \begin{theorem}\label{thm:2nd-order-i-group-from-affine-connection}
    Let $M$ be a formal manifold and $\lambda$ an affine connection on $M$; then $\lambda$ induces an i-group structure on $\frakM_2(P)$ for each $P\in M$ such that $\lambda(P,Q,R)=QR$ for $Q,R\in\frakM(P)$.
  \end{theorem}
  \begin{proof}
    The proof combines the observation about torsion with the proof strategies of theorem~\ref{thm:lie-algebra-from-i-algebra}, theorem~\ref{thm:formal-2nd-order-i-groups} and proposition~\ref{prop:symmetric_affine_connection}. 
    
    \begin{enumerate}[(i)]
      \item Firstly, let $\lambda^{op}$ denote the conjugated affine connection of $\lambda$, i.e. $\lambda^{op}(P,Q,R) = \lambda(P,R,Q)$. The symmetrisation 
      $$
        \bar{\lambda} = \frac{1}{2}\lambda + \frac{1}{2}\lambda^{op}
      $$  
      yields a symmetric affine connection and thus a second-order i-affine structure on $M_2$ by proposition~\ref{prop:symmetric_affine_connection}.
      
      \item We use the torsion $\tau$ of $\lambda$ to construct an alternating bilinear map $[-,-]_P$ on $T_PM$ characterised by
      $$
        [t_1,t_2]_P(d_1d_2) = -\tau_P\bigl(t_1(d_1),t_2(d_2)\bigr),\qquad d_1,d_2\in D
      $$
      as in theorem~\ref{thm:lie-algebra-from-i-algebra}.
      
      The representation of torsion in a chart \cite[eq.~2.3.13]{Kock:Synthetic_Geometry_Manifolds} 
      \begin{equation}\label{eq:torsion-in-chart}
        \tau_{P}(Q,R)= P - \bigl(\Gamma_P[Q-P,R-P]-\Gamma_P[R-P,Q-P]\bigr),
      \end{equation}
      where $\Gamma_P$ denotes the connection symbol of $\lambda$, shows that $\tau_P(Q,R) \in\frakM(P)$ for $Q,R\in \frakM(P)$, and hence $-\tau_P(Q,R)\in\frakM(P)$. We get an alternating map
      $$
        \frakM(P)\times \frakM(P) \to \frakM(P),\qquad (P,Q)\mapsto -\tau_P(P,Q)
      $$
      that lifts to an alternating map
      $$
        [-,-]_P: D(T_PM)\times D(T_PM) \to D(T_PM)
      $$
      via the first-order $\exp$-$\log$ bijection at $P$. Since $\tau_P(P,Q)=\tau_P(Q,P) = P$ we have $[0,t]_P=[t,0]_P=0$. By the Kock-Lawvere axiom~(I) $[-,-]_P$ has thus a unique extension to an alternating bilinear map on $T_PM$. By considering the bracket in a chart we find
      $$
        [t_1,t_2]_P(d_1d_2) = -\tau_P\bigl(t_1(d_1),t_2(d_2)\bigr),\qquad d_1,d_2\in D,
      $$
      which characterises $[t_1,t_2]_P$ due to the Kock-Lawvere axiom~(I) as explained in the proof of theorem~\ref{thm:lie-algebra-from-i-algebra}.

      \item By theorem~\ref{thm:formal-2nd-order-i-groups} the operations
      $$
        (t_1,t_2)\mapsto t_1+t_2+\frac{1}{2}[t_1,t_2]_P
      $$
      and
      $$ 
        t\mapsto -t
      $$
      make $D_2(T_PM)$ into an i-group. We transport this i-group structure to $\frakM_2(P)$ using the second-order $\log$-$\exp$ bijection (\ref{eq:log_exp}), i.e.
      $$
        QR := \exp_P\bigl(\log_P(Q)\log_P(R)\bigr) 
      $$
      and
      $$
        Q^{-1} := \exp_P\bigl(-\log_P(Q)\bigr)
      $$
      for $\la Q,R\ra$ in $\frakM_2(P)$. This makes $\frakM_2(P)$ into an i-group. 
      
      \item Let $Q,R\in\frakM(P)$. By lemma~\ref{lem:embedding-into-2nd-order} we have $\la Q,R\ra$ in $\frakM_2(P)$. Since $\exp_P$ is i-linear for the second-order i-linear structure we find
      $$
        Q^{-1} = \exp_P\bigl(-\log_P(Q)\bigr) = -Q
      $$
      and
      \begin{align*}
        QR &= \exp_P\Bigl(\log_P(Q)+\log_P(R)+\frac{1}{2}\bigl[\log_P(Q),\log_P(R)\bigr]_P\Bigr)\\
        &= Q + R + \frac{1}{2}\exp_P\Bigl(\bigl[\log_P(Q),\log_P(R)\bigr]_P\Bigr) \\
        &= Q + R - \frac{1}{2}\tau_P(Q,R),
      \end{align*}
      where the last equation follows from the construction of the bracket $[-,-]_P$ on $T_PM$. The last i-linear combination in $\frakM_2(e)$ can be re-written as the i-affine combination in $M_2$
      $$ 
        QR = Q + R - P - \frac{1}{2}\bigl(\tau_P(Q,R) - P\bigr)
      $$
      Using the presentation of torsion in a chart~(\ref{eq:torsion-in-chart}) a calculation like the one in the proof of theorem~\ref{thm:2nd-order-operations}~(2) shows that 
      \begin{align*}
        QR &= Q + R - P + \bar{\Gamma}_P[Q-P,R-P] + \frac{1}{2}\bigl(\Gamma_P[Q-P,R-P] - \Gamma_P[R-P,Q-P]\bigr)\\
        &= Q + R - P + \Gamma_P[Q-P,R-P] \\
        &= \lambda(P,Q,R),
      \end{align*}
      where $\bar{\Gamma}_P$ is the connection symbol of $\bar{\lambda}$, and thus the symmetrisation of $\Gamma_P$. This shows that the group operation agrees with $\lambda$ as claimed.
    \end{enumerate}
      
  \end{proof}

  Comparing theorem~\ref{thm:2nd-order-i-group-from-affine-connection} with proposition~\ref{prop:symmetric_affine_connection} we notice the absence of a notion of algebraic structure that does not rely on base points. Indeed, proposition~\ref{prop:symmetric_affine_connection} implies that every symmetric affine connection on $M$ is equivalent to the existence of a second-order i-linear structure on every $\frakM_2(P)$ and $P\in M$. It is only when passing from i-linear to i-affine combinations that the operations become base-point independent, and we speak of an i-affine structure. 
  
  A direct generalisation to the non-abelian case would be the following operations of \emph{non-abelian i-affine combinations} on $M_2$:
  $$
    \calA(n)\times M_2\la n\ra,\quad (\mu,\la P_1,\ldots,P_n\ra) \mapsto \exp_{P_1}\Bigl(\prod_{j=1}^n \mu_j\log_{P_1}(P_j)\Bigr)
  $$
  It turns out, however, that torsion forms an obstruction to their base-point independence.  
  \begin{proposition}
    Let $\la P,Q,P_1,\ldots,P_n\ra\in M_2\la n+2\ra$ and $\mu \in \calA(n)$, then
    \begin{equation*}
      \begin{split}
        \exp_{P}\Bigl(\prod_{j=1}^n \mu_j\log_{P_1}(P_j)\Bigr) = &\exp_{Q}\Bigl(\prod_{j=1}^n \mu_j\log_{P_2}(P_j)\Bigr)  \\ &+ \frac{1}{2}\sum_{j=1}^n\sum_{k<j}\mu_k\mu_j\bigl(\tau_Q(P,P_k) +\tau_Q(P_j,P) -2Q\bigr)
      \end{split}
    \end{equation*}
  \end{proposition}
  \begin{proof}
    Firstly, we represent and simplify the $n$-ary operation using the second-order i-affine structure induced by $\bar{\lambda}$ on $M_2$ to
    $$
      \exp_{P}\Bigl(\prod_{j=1}^n \mu_j\log_{P_1}(P_j)\Bigr) = \sum_{j=1}^n \mu_j P_j - \frac{1}{2}\sum_{j=1}^n\sum_{k<j}\mu_k\mu_j\bigl(\tau_P(P_k,P_j)-P\bigr)
    $$
    respectively,
    $$
      \exp_{P}\Bigl(\prod_{j=1}^n \mu_j\log_{P_1}(P_j)\Bigr) = \sum_{j=1}^n \mu_j P_j - \frac{1}{2}\sum_{j=1}^n\sum_{k<j}\mu_k\mu_j\bigl(\tau_Q(P_k,P_j)-Q\bigr) 
    $$
    (Here we have identified $\tau$ with its unique extension to $M_2\la 3\ra$ as an alternating i-biaffine map in the second and third argument, guaranteed by theorem~\ref{thm:2nd-order-i-group-from-affine-connection} and Kock-Lawvere axiom~(III).)
    
    We shall now show 
    $$
      \tau_P(R,S) - P= \tau_Q(R,S) + \tau_Q(P,R) + \tau_Q(S,P)- 3Q
    $$
    for $\la P,Q,R,S\ra$ in $M_2$. We work in a chart. Applying (\ref{eq:torsion-in-chart}) yields
    $$
      \tau_P(R,S) - P = \Gamma_P[S-P,R-P] - \Gamma_P[R-P,S-P]
    $$
    Taylor expansion of $\Gamma_P = \Gamma_{P+(Q-P)}$ and bearing in mind that third and higher-order terms vanish results in
    \begin{align*}
      \Gamma_Q[S-P,R-P] &= \Gamma_P[S-P,R-P] +\partial \Gamma_P[S-P,R-P][Q - P] \\
      &\qquad\qquad + \frac{1}{2}\partial^2 \Gamma_P[S-P,R-P][Q - P]^2\\
      &= \Gamma_P[S-P,R-P].
    \end{align*}
    Finally, we apply the bilinearity of $\Gamma_Q$: 
    \begin{align*}
      \tau_P(R,S) - P &= \Gamma_P[S-P,R-P] - \Gamma_P[R-P,S-P]\\
      &= \Gamma_Q[S-Q + Q - P,R - Q + Q - P] \\
      &\qquad\qquad - \Gamma_P[R - Q + Q - P,S - Q + Q - P]\\
      &= \Gamma_Q[S-Q,R-Q] - \Gamma_Q[R-Q,S-Q] \\
      &\qquad\qquad - \Gamma_Q[S-Q,P-Q] - \Gamma_Q[P-Q,R-Q] \\
      &\qquad\qquad + \Gamma_Q[R-Q,P-Q] + \Gamma_Q[P-Q,S-Q]  \\
      &= \tau_Q(R,S) - Q + \tau_Q(P,R) - Q + \tau_Q(S,P)- Q.
    \end{align*}
    Note that the last equality holds in $M_2$, which follows from the fact that \\ $\la Q,\tau_Q(R,S), \tau_Q(P,R), \tau_Q(S,P)\ra$ in $M_1\into M_2$ (cf. proposition~\ref{prop:lie-bracket-of-points}).
  \end{proof}

%%%%%%%%%%%%%%%%%%%%%%%%%%%%%%%%%%%%%%%%%%%%%%%%%%%%%%%%%%%%%%%%%%%%%%%%%%%%%%%
\section{Conclusion}

Going beyond infinitesimally affine spaces and vector spaces we have introduced the notion of an infinitesimal group as an infinitesimal model of the algebraic theory of groups. We have shown that for a Lie group the first- and second-order monads of the neutral element are i-subgroups of the Lie group $G$; but that this is not the case for the monad induced by the nil-square i-structure. This shows that the naturally defined first-order i-structure becomes relevant in the case of non-commutative infinitesimal algebra, where the nil-square i-structure is too big.

With the Lie bracket of points we were able to re-construct the Lie bracket of the Lie algebra of a Lie group. Since the second-order i-group operation is i-biaffine we have seen that the Lie bracket can always be understood as a commutator familiar from matrix Lie groups, even in the absence of an enveloping algebra. Moreover, we were able to characterise second-order i-group structures on KL vector spaces as deformations of vector addition by bilinear maps.  

Moving from the infinitesimal to the global picture, with the relationship between torsion and the Lie bracket of points we were able to generalise the construction of second-order i-groups from Lie groups to manifolds with non-symmetric affine connections. Since this reduces to the second-order i-affine structure in the case of vanishing torsion, it also provides a generalisation of the correspondence between symmetric affine connections and second-order i-affine structures on a manifold established in \cite{Bar:second_order_affine_structures}. However, unlike in the symmetric case we are currently lacking the corresponding algebraic theory that would induce a (non-symmetric) affine connection, which remains an open problem.

An important example we have not discussed in this paper is the diffeomorphism group of a formal manifold. Although it is not a formal manifold, in general, the endomorphism monoid is an i-affine space for the nil-square i-structure \cite[chap.~3.3.2]{Bar:thesis} and a microlinear space \cite[chap.~2.3]{Lavendhomme:Basic_Concepts_of_SDG}. We shall discuss this in a future paper. 

Last but not least, we have seen how infinitesimal algebra allows to simplify and clarify many of the basic differential-geometric constructions on Lie groups and formal manifolds. It is re-assuring to see how familiar algebraic structures simply re-appear on the infinitesimal level and interplay with the global structure like in the case of Lie groups, for example. On the infinitesimal level, many arguments in the chart are not for coordinatisation, but for doing the multilinear algebra. This formal multilinear algebra simplifies considerably when it can be translated into infinitesimal algebra of points. We believe that by developing infinitesimal algebra further, many calculations involving infinitesimals will eventually be able to be done on the manifold itself, not requiring a chart. 

\bibliography{references_with_links}

\end{document}